# AUTOMORPHISMS OF THE CYCLE PREFIX DIGRAPH

William Y. C. Chen
Center for Combinatorics
The Key Laboratory of Pure Mathematics and Combinatorics of Ministry of Education
Nankai University, Tianjin 300071, P. R. China

Vance Faber
Big Pine Key, FL

Bingqing Li
Department of Risk Management and Insurance
Nankai University, Tianjin 300071, P. R. China

**AMS Classification**: 11B83, 05C05

**Keywords**: Cycle prefix digraph, automorphism groups

**Note**. This is a version of an unpublished report from 2001. We place it here to give it wider access.

**Abstract.** Cycle prefix digraphs have been proposed as an efficient model of symmetric interconnection networks for parallel architecture. It has been discovered that the cycle prefix networks have many attractive communication properties. In this paper, we determine the automorphism group of the cycle prefix digraphs. We show that the automorphism group of a cycle prefix digraph is isomorphic to the symmetric group on its underlying alphabet. Our method can be applied to other classes of graphs built on alphabets including the hypercube, the Kautz graph,and the de Bruijn graph.

Keywords: Cycle prefix digraph, automorphism group

AMS Subject Classification: 05C25, 68M10

# 1 Introduction

The cycle prefix digraphs $\Gamma_\Delta(D)$,or the cycle prefix networks, have been recently proposed as model of symmetric interconnection networks with a large number of vertices and small diameter [5]. It has been discovered that such networks possess a number of attractive properties for efficient communication in parallel computing. The main advantage of using a symmetric network is that each processor sees the network communication the same as any other processor in the network [1]. This implies that all the processors can essentially execute the same communication program. For this reason, a symmetric network is sometimes called a homogeneous network. Beyond the vertex symmetry property of a digraph, the automorphism group actually provides a complete picture of the symmetries of the digraph. An important characterization of symmetric digraphs is due to Sabidussi which states that a digraph is vertex symmetric if and only if it can be represented by a Cayley coset digraph [9]. This characterization leads to the discovery of the cycle prefix digraph $\Gamma_\Delta(D)$ $(\Delta \geq D)$, and the reader is referred to [5] for the Cayley coset digraph definition of $\Gamma_\Delta(D)$. However, in the present paper, we shall use the representation of $\Gamma_\Delta(D)$ in terms of sequences, analogous to the representations of hypercubes, the de Bruijn graphs and the Kautz graphs.

Throughout this paper, we assume that $\Delta \geq D$. The sequence definition of $\Gamma_\Delta(D)$ goes as follows. A vertex of $\Gamma_\Delta(D)$ is represented by a partial permutation $x_1x_2\ldots x_D$ on $1,2,\ldots,\Delta+1$, which is a sequence of distinct elements. A partial permutation of $D$ elements is also called a $D$-permutation. Since the adjacency in a digraph could be ambiguous if one does not specify the direction, if $(u,v)$ is an arc in a digraph, we shall say that $u$ is adjacent to $v$, whereas $v$ is next to $u$. The adjacency relation of $\Gamma_\Delta(D)$ can be described as follows:

$$x_1x_2\ldots x_D \Rightarrow \begin{cases} x_kx_1\ldots x_{k-1}\ldots x_D, & \text{if} \quad 2\leq k\leq D; \\ mx_1x_2\ldots x_{D-1}, & \text{if} \quad m\neq x_1,x_2,\ldots,x_{D.} \end{cases}$$

We say that $x_kx_1\ldots x_{k-1}\ldots x_D$ is next to $x_1x_2\ldots x_D$ via a rotation operation for $2\leq k\leq D$, and $mx_1x_2\ldots x_{D-1}$ is said to be next $x_1x_2\ldots x_D$ via a shift operation for $m\neq x_1,x_2,\ldots,x_{D.}$. In particular, if $k=D$, the vertex $x_Dx_1x_2\ldots x_{D-1}$ is said to be next to $x_1x_2\ldots x_D$ via a full rotation. Moreover, we shall use the following notation:

$$R_k(x_1x_2\ldots x_D)=x_kx_1\ldots x_{k-1}\ldots x_D \qquad \text{for } 2\leq k\leq D;$$

$$R_m(x_1x_2\ldots x_D)=mx_1x_2\ldots x_{D-1} \qquad \text{for } m\neq x_1,x_2,\ldots,x_{D.}$$

In this paper, we will discuss the general case $\Gamma_\Delta(D,-r)$, $r\geq 0$, in which the adjacency relation is defined as:

$$x_1x_2\ldots x_D \Rightarrow \begin{cases} x_k x_1 \ldots x_{k-1} \ldots x_D, & if \quad r+2 \le k \le D; \\ mx_1x_2\ldots x_{D-1}, & if \quad m \ne x_1, x_2, \ldots, x_{D.} \end{cases}$$

Note that $\Gamma_\Delta(D)$ is a special case of $\Gamma_\Delta(D,-r))$ for $r=0$.

From the above sequence representation, one sees that $\Gamma_\Delta(D,-r)$ is vertex symmetric [5]: Given any two vertices $X = x_1x_2\ldots x_D$ and $Y = y_1y_2\ldots y_D$ ,one can always find a permutation $\pi$ on $1,2,\ldots,\Delta+1$ such that $\pi(x_i) = y_i$. Clearly, such a permutation $\pi$ induces an automorphism of $\Gamma_\Delta(D,-r)$ because it simply relabels the underlying alphabet of $\Gamma_\Delta(D,-r)$. It is easy to see that every vertex in $\Gamma_\Delta(D,-r)$ has both indegree and outdegree $\Delta - r$ , and the digraph has $(\Delta+1)_D = (\Delta+1)\Delta\cdots(\Delta - D + 2)$ vertices. Moreover, $\Gamma_\Delta(D,-r)$ has diameter $D + r$ [3] and many other nice properties including Hamiltonicity, high connectivity, small wide diameter, and regular reachability [2,3,4,5,7,6,8].

The main concern of this paper is with the automorphism group of $\Gamma_\Delta(D,-r)$. We shall show that the automorphism group of $\Gamma_\Delta(D,-r)$ is isomorphic to the symmetric group on its underlying set $1,2,\ldots,\Delta+1$. Our approach to determine the automorphism group of the cycle prefix digraphs may be employed to characterize the automorphism group of other classes of networks,such as the hypercube, the Kautz graph and the de Bruijn graph.

## 2 The Automorphism Group of $\Gamma_\Delta(D,-r)$

Given a digraph $G = (V,E)$, an automorphism of $G$ is a permutation $\pi$ on $V$ such that $(u,v) \in E$ if and only if $(\pi(u),\pi(v)) \in E$ . All the automorphisms of $G$ form the automorphism group of $G$. We say that $G$ is symmetric or vertex transitive if every vertex of $G$ looks the same - strictly speaking, for any two vertices $u$ and $v$ there exists an automorphism $\pi$ of $G$ such that $\pi(u)=v$. The main result of this paper is to determine the automorphism group of the cycle prefix digraph $\Gamma_\Delta(D,-r)$. Namely, we will prove the following:

**Theorem 2.1.** *For $\Delta \geq D \geq r+2$, the automorphism group of $\Gamma_\Delta(D,-r)$ is isomorphic to the symmetric group $S_{\Delta+1}$ on the underlying set $\{1,2,\ldots,\Delta+1\}$.*

We establish a few lemmas.

**Lemma 2.2.** *Let $\alpha$ be an automorphism of $\Gamma_\Delta(D,-r)$ and $u$, $v$ two vertices in $\Gamma_\Delta(D,-r)$. If $v$ is next to $u$ via a shift operation then $\alpha(v)$ is next to $\alpha(v)$ via a shift operation as well.*

*Proof:* Since $\Gamma_\Delta(D,-r)$ is vertex symmetric, without loss of generality we may assume that $u=12\ldots D$ and $v=i12\ldots(D-1)$, where $i$ is an integer satisfying $D<i\leq\Delta+1$. Let $\alpha(12\ldots D)=w_1w_2\ldots w_D$ and let $\alpha(i12\ldots(D-1))=y_1y_2\ldots y_D$. Since $i12\ldots(D-1)$ is next to $12\ldots D$ in $\Gamma_\Delta(D)$ and $\alpha$ is an automorphism, $y_1y_2\ldots y_D$ is next $w_1w_2\ldots w_D$ in $\Gamma_\Delta(D)$. Note that the distance from $i12\ldots(D-1)$ to $12\ldots D$ is at least $D$ because $D$ does not appear in $i12\ldots(D-1)$ and either a rotation or shift operation can move an element in a sequence at most one step forward (from left to right). Assume that $y_1y_2\ldots y_D$ is next to $w_1w_2\ldots w_D$ via a rotation operation, say, $R_k$ for some $2\leq k<D$, that is, $y_1y_2\ldots y_D=R_k(w_1w_2\ldots w_D)$. Then it is clear that $R_k^k(w_1w_2\ldots w_D)=w_1w_2\ldots w_D$. It follows that $R_k^{k-1}(y_1y_2\ldots y_D)=w_1w_2\ldots w_D$. Therefore, the distance from $y_1y_2\ldots y_D$ to $w_1w_2\ldots w_D$ is at most $k-1$. Since an automorphism preserves distance, $\alpha$ *be an automorphism* must be next to $w_1w_2\ldots w_D$ via a shift operation. □

Since the adjacency of two vertices in $\Gamma_\Delta(D,-r)$ is realized either by a rotation or shift operation, Lemma 2.2 is equivalent to the statement that if $v$ is next to $u$ via a rotation operation, then $\alpha(v)$ is next to $\alpha(v)$ via a rotation operation as well.

**Lemma 2.3.** *Let $\alpha$ be an automorphism and $x_1x_2\ldots x_D$, $y_1y_2\ldots y_D$ be two vertices of $\Gamma_\Delta(D,-r)$ such that $y_1y_2\ldots y_D = \alpha(x_1x_2\ldots x_D)$. Let $S_m$ be a shift operation with $m \in \{1,2,\cdots,\Delta+1\}\setminus\{x_1,x_2,\cdots,x_D\}$. If $x_i$, $1 \le i \le D$, appears in $S_m(x_1x_2\ldots x_D)$, then $y_i$ must appear in $\alpha(S_m(x_1x_2\ldots x_D))$ at the same position as $x_i$ in $S_m(x_1x_2\ldots x_D)$.*

*Proof.* Since $y_1y_2\ldots y_D = \alpha(x_1x_2\ldots x_D)$, by Lemma 2.2 we have $\alpha(S_m(x_1x_2\ldots x_D)) = ty_1y_2\ldots y_{D-1}$ for some $t$. Clearly, $x_i$ and $y_i$ are at the same position for $1 \le i \le D$. □

**Lemma 2.4.** *Let $\alpha$ be an automorphism and $x_1x_2\ldots x_D$, $y_1y_2\ldots y_D$ be two vertices of $\Gamma_\Delta(D,-r)$ such that $y_1y_2\ldots y_D = \alpha(x_1x_2\ldots x_D)$.. Let $R_k$, $r+2 \le k \le D$, be a rotation operation. If $x_i$, $1 \le i \le D$, appears in $R_k(x_1x_2\ldots x_D)$, then $y_i$ must appear in $\alpha(R_k(x_1x_2\ldots x_D))$ at the same position as $x_i$ in $R_k(x_1x_2\ldots x_D)$.*

*Proof.* Let

$$w_1w_2\ldots w_D = R_k(x_1x_2\ldots x_D) = x_kx_1\ldots x_{k-1}x_{k+1}\ldots x_D$$

(2.1)

From Lemma 2.2, it follows that $\alpha(w_1w_2\ldots w_D)$ is next to $\alpha(x_1x_2\ldots x_D)$ via a rotation operation,say,

$$\alpha(w_1w_2\ldots w_D)=R_j(\alpha(x_1x_2\ldots x_D))=R_j(y_1y_2\ldots y_D),$$

where $r+2\le j\le D$. Since $R_k(x_1x_2\ldots x_D)=x_kx_1\ldots x_{k-1}x_{k+1}\ldots x_D$, it takes at least $k-1$ steps to move $x_k$ from the first position to the $k$ th position. Hence the distance from $R_k(x_1x_2\ldots x_D)$ to $x_1x_2\ldots x_D$ is at least $k-1$. On the other hand, we have $R_k^{k-1}(R_k(x_1x_2\ldots x_D))=x_1x_2\ldots x_D$. Thus, the distance from $R_k(x_1x_2\ldots x_D)$ to $x_1x_2\ldots x_D$ is exactly $k-1$. Similarly,the distance from $\alpha(w_1w_2\ldots w_D)=R_j(y_1y_2\ldots y_D)$ is $j-1$. Since an automorphism preserves distance,we have $k=j$. Therefore, we have

$$\alpha(w_1w_2\ldots w_D)=R_k(y_1y_2\ldots y_D)=y_ky_1\ldots y_{k-1}y_{k+1}\ldots y_D \tag{2.2}$$

We now reach the claim by comparing (2.1) with (2.2).

We may illustrate Lemmas 2.3 and 2.4 by the following diagram:

$$\begin{array}{ccc} x_1x_2\ldots x_D & \overset{S\ or\ R}{\rightarrow} & \ldots x_i\ldots \\ \downarrow\alpha & & \downarrow\alpha \\ y_1y_2\ldots y_D & \overset{S\ or\ R}{\rightarrow} & \ldots y_i\ldots \end{array} \tag{2.3}$$

**Lemma 2.5.** *Let* $\alpha$ *be an automorphism of* $\Gamma_\Delta(D,-r)$. *Then there exists a permutation* $w = w_1 w_2 \ldots w_{\Delta+1}$ *on the underlying set* $\{1,2,\cdots,\Delta+1\}$ *of* $\Gamma_\Delta(D,-r)$ *such that*

$$\alpha(12\ldots D) = w_1 w_2 \ldots w_D, \tag{2.4}$$

$$\alpha(i12\ldots D-1) = w_i w_1 w_2 \ldots w_{D-1} \tag{2.5}$$

*for* $i = D+1, D+2, \ldots, \Delta+1$.

*Proof.* We choose $w_1 w_2 \ldots w_D$ as the image of $12\ldots D$ under the action of $\alpha$. For $D \le i \le \Delta+1$, since $i12\ldots(D-1)$ is next to $12\ldots D$ via a shift operation, then $\alpha(i12\ldots(D-1))$ is next to $\alpha(12\ldots D) = w_1 w_2 \ldots w_D$ via a shift operation as signified in (2.5). Since $\alpha$ is a one-to-one map on the vertices of $\Gamma_\Delta(D,-r)$, the elements $w_{D+1}, \ldots, w_{\Delta+1}$ are distinct to each other as well as to $w_1, w_2, \ldots, w_D$. It follows that the above chosen sequence $w$ is a permutation on $\{1,2,\cdots,\Delta+1\}$. □

We call the permutation $w$ constructed from the automorphism $\alpha$ in the above lemma the *derived permutation* of $\alpha$. Our goal is to show that the automorphism $\alpha$ is uniquely determined by its derived permutation. For convenience, we also write $w_i$ as $w(i)$. Let us illustrate the above lemma by the following diagram

$$\begin{array}{ccc} 12\ldots D & \overset{S_i}{\rightarrow} & i12\ldots(D-1) \\ \downarrow \alpha & & \downarrow \alpha \\ w_1 w_2 \ldots w_D & \overset{S_{w_i}}{\rightarrow} & w_i w_1 w_2 \ldots w_{D-1} \end{array} \tag{2.6}$$

**Definition 2.6.** *Let* $\alpha$ *be an automorphism of* $\Gamma_\Delta(D,-r)$ *and* $w = w_1 w_2 \ldots w_{\Delta+1}$ *a permutation on* $\{1,2,\cdots,\Delta+1\}$. *Given a vertex* $x_1 x_2 \ldots x_D$ *in* $\Gamma_\Delta(D,-r)$, *we say that* $\alpha$ *and* $w$ *are compatible for* $x_1 x_2 \ldots x_D$ *if* $\alpha$ *and* $w$ *have the same image of* $x_1 x_2 \ldots x_D$, *namely,*

$$\alpha(x_1 x_2 \ldots x_D) = w(x_1) w(x_2) \ldots w(x_D) \ . \tag{2.7}$$

**Lemma 2.7.** *Let α be an automorphism of* $\Gamma_\Delta(D,-r)$ *and* $w = w_1, w_2, \ldots, w_{\Delta+1}$ *be its derived permutation. Suppose* $\alpha$ *and* $w$ *are compatible for a vertex* $x_1 x_2 \ldots x_D$. *If* $y_1 y_2 \ldots y_D$ *is next to* $x_1 x_2 \ldots x_D$ *via a rotation operation, then* $\alpha$ *and* $w$ *are compatible for* $y_1 y_2 \ldots y_D$ *as well.*

*Proof.* Suppose $\alpha(x_1 x_2 \ldots x_D) = w(x_1) w(x_2) \ldots w(x_D)$ and $y_1 y_2 \ldots y_D = R_k(x_1 x_2 \ldots x_D)$, $r+2 \le j \le D$. By Lemma 2.4, we have the following diagram:

$$\begin{array}{ccc} x_1 x_2 \ldots x_D & \xrightarrow{R_k} & R_k(x_1 x_2 \ldots x_D) \\ \downarrow \alpha & & \downarrow \alpha \\ w(x_1) w(x_2) \ldots w(x_D) & \xrightarrow{R_k} & R_k(w(x_1) w(x_2) \ldots w(x_D)) \end{array} \tag{2.8}$$

That is, $\alpha(y_1 y_2 \ldots y_D) = w(y_1) w(y_2) \ldots w(y_D)$ as desired. □

**Lemma 2.8.** *Let α be an automorphism of* $\Gamma_\Delta(D,-r)$, $\Delta \ge D \ge r+2$, *and* $w = w_1 w_2 \ldots w_{\Delta+1}$ *be its derived permutation. Suppose* $\alpha$ *and* $w$ *are compatible* *for a vertex* $x_1x_2\ldots x_D$*. If* $y_1y_2\ldots y_D$ *is next to* $x_1x_2\ldots x_D$ *via a shift operation, then* $\alpha$ *and* $w$ *are compatible for* $y_1y_2\ldots y_D$ *as well.*

*Proof.* From the given condition, we may write $y_1y_2\ldots y_D$ as $ix_1x_2\ldots x_{D-1}$, $i \notin \{x_1, x_2, \cdots, x_{D.}\}$. Suppose $\alpha(x_1x_2\ldots x_D) = w(x_1)w(x_2)\ldots w(x_D)$. By Lemma 2.3, we have the following diagram:

$$\begin{array}{ccc}
x_1x_2\ldots x_D & \overset{S_i}{\rightarrow} & ix_1x_2\ldots x_{D-1} \\
\downarrow\alpha & & \downarrow\alpha \\
w(x_1)w(x_2)\ldots w(x_D) & \overset{S_j}{\rightarrow} & jw(x_1)w(x_2)\ldots w(x_{D-1})
\end{array}$$

.

It remains to show that $j = w(i)$. For convenience, we say that a path contains $i$ if $i$ is contained in every vertex in the path. There are two cases:

Case 1. $i > D$. By Lemma 2.3, 2.4 and 2.5, if we can find a path from $ix_1x_2\ldots x_{D-1}$ to $i12\ldots(D-1)$ that contains $i$, then we can get the following diagram:

$$\begin{array}{ccccc}
x_1x_2\ldots x_D & \overset{S_i}{\rightarrow} & ix_1x_2\ldots x_{D-1} & \overset{\ldots}{\rightarrow} & i12\ldots(D-1) \\
\downarrow\alpha & & \downarrow\alpha & & \downarrow\alpha \\
w(x_1)\ldots w(x_D) & \overset{S_j}{\rightarrow} & jw(x_1)\ldots w(x_{D-1}) & \overset{\ldots}{\rightarrow} & w(i)w(1)\ldots w(D-1)
\end{array}$$

From this diagram, we obtain $j = w(i)$. Now we show that there is a path from $ix_1x_2\ldots x_{D-1}$ to $i12\ldots(D-1)$ that contains $i$. In fact, if $D-1 \in \{x_1, x_2, \cdots, x_{D-1}\}$, say $x_j = D-1$, then $(D-1)ix_1\ldots x_{j-1}x_{j+1}\ldots x_{D-1}$ is next to $ix_1x_2\ldots x_{D-1}$ via a rotation operation. If $D-1 \notin \{x_1, x_2, \cdots, x_{D-1}\}$, then $(D-1)ix_1x_2\ldots x_{D-2}$ is next to $ix_1x_2\ldots x_{D-1}$ via a shift operation. Step by step,finally we can get the following path:

$$\begin{aligned} &ix_1x_2\ldots x_{D-1} \to (D-1)i\cdots \to (D-2)(D-1)i\cdots \to \cdots \\ &\qquad 12\ldots(D-1)i \to i12\ldots(D-1) \end{aligned} \tag{2.9}$$

which is a path $ix_1x_2\ldots x_{D-1}$ to $i12\ldots(D-1)$ that contains $i$.

Case 2. $i \le D$. Analogous to the above case, if we can find a path from $ix_1x_2\ldots x_{D-1}$ to $i12\ldots(D-1)$ that contains $i$, then we have the following Diagram

$$\begin{array}{ccccc} x_1x_2\ldots x_D & \overset{S_i}{\to} & ix_1x_2\ldots x_{D-1} & \overset{\cdots}{\to} & 12\ldots i\ldots D \\ \downarrow\alpha & & \downarrow\alpha & & \downarrow\alpha \\ w(x_1)\ldots w(x_D) & \overset{S_j}{\to} & jw(x_1)\ldots w(x_{D-1}) & \overset{\cdots}{\to} & w(1)\ldots w(i)\ldots w(D) \end{array}$$

Similar to Case 1, we may find a path from $ix_1x_2\ldots x_{D-1}$ to $12\ldots i\ldots D$ that contains $i$.

Combining the above two cases, Lemma 2.8 is proved. □

We are now ready to prove Theorem 2.1.

*Proof of Theorem 2.1*. Since each permutation $w$ on $\{1,2,\cdots,\Delta+1\}$ induces an automorphism $\alpha$ of $\Gamma_\Delta(D,-r)$ with $w$ being the derived permutation of $\alpha$, we only need to show that any automorphism $\alpha$ of $\Gamma_\Delta(D,-r)$ is uniquely determined by its derived permutation. Let $\alpha$ be an automorphism of $\Gamma_\Delta(D,-r)$ with derived permutation $w = w_1, w_2, \ldots, w_{\Delta+1}$. We aim to show that $\alpha$ and $w$ are compatible for any vertex $x_1 x_2 \ldots x_D$ of $\Gamma_\Delta(D,-r)$, that is,

$$\alpha(x_1 x_2 \ldots x_D) = w(x_1) w(x_2) \ldots w(x_D)\,. \tag{2.10}$$

Our strategy for proving (2.10) is as follows. By definition, it holds for the initial vertex $12\ldots D$. Assuming that it holds for a a vertex $x_1 x_2 \ldots x_D$, if we can show that it also holds for any vertex $y_1 y_2 \ldots y_D$ next to $x_1 x_2 \ldots x_D$, then (2.10) must be true for any vertex because $\Gamma_\Delta(D,-r)$ is strongly connected.

Let us assume that $x_1 x_2 \ldots x_D$ satisfies (2.10). Suppose $y_1 y_2 \ldots y_D$ is a vertex next to $x_1 x_2 \ldots x_D$. If $y_1 y_2 \ldots y_D$ is next to $x_1 x_2 \ldots x_D$ via a rotation operation, from Lemma 2.7 it follows that (2.10) holds for $y_1 y_2 \ldots y_D$. When $y_1 y_2 \ldots y_D$ is next to $x_1 x_2 \ldots x_D$ via a shift operation, from Lemma 2.8 it follows that (2.10) also holds for $y_1 y_2 \ldots y_D$. Noting that different derived permutations lead to different automorphisms, this completes the proof.

Let us examine some special cases of the above theorem. Let $n = \Delta + 1$. When $\Delta = D$, $\Gamma_\Delta(D)$ turns out to be the Cayley graph on the symmetric group $S_n$ with generators $(12)$, $(123)$, …, $(12\ldots n)$. In this case, it is much easier to see that the automorphism group is $S_n$ because no shift operation is involved.

The simplest nontrivial case is for $D = 2$. Note that $\Gamma_\Delta(D)$ is the Kautz graph $K(\Delta,2)$. It is known that the automorphism group of the Kautz graph $K(\Delta, D)$ is isomorphic to the symmetric group on the alphabet $\{1,2,\cdots,\Delta+1\}$ (see [10]).

**Acknowledgment.** This work was partially supported by the 973 project sponsored by the Ministry of Science of Technology of China, and the National Science Foundation of China. The authors would like to thank Dr. Q. L. Li for valuable comments.